\documentclass[11pt]{article}
\usepackage{amsmath,amssymb,mathrsfs}

\newcommand{\R}{\mathbb{R}}
\newcommand{\Z}{\mathbb{Z}}
\newcommand{\p}{\mathbb{P}}
\newcommand{\E}{\mathbb{E}}
\newcommand{\D}{\mathcal{D}_{\epsilon}(\omega)}

\newcommand{\B}{\mathcal{B}_{\epsilon}(\omega)}
\newcommand{\A}{\mathcal{A}}
\newcommand{\T}{\mathcal{T}}
\renewcommand{\index}{i \mspace{-3mu} \imath}
\newcommand{\Ocl}{\mathcal{O}_{\epsilon}(\omega)}
\newcommand{\den}{\mathcal{D}_{\epsilon}(\omega_1,\omega_2)}
\newcommand{\clden}{\overline{\mathcal{D}}_{\epsilon}(\omega_1,\omega_2)}
\newcommand{\bad}{\mathcal{B}_{\epsilon}(\omega_1,\omega_2)}
\newcommand{\Oc}{\mathcal{O}_{\epsilon}(\omega_1,\omega_2)}
\newcommand{\set}[1]{\left\{#1\right\}}
\newcommand{\braket}[1]{\left\langle#1\right\rangle}
\usepackage{theorem}
\theorembodyfont{\upshape}
\newtheorem{thm}{Theorem}
\newtheorem{lem}{Lemma}
\newtheorem{prop}{Proposition}
\newtheorem{cor}{Corollary}
\newtheorem{Def}{Definition}

\begin{document}

\title{Asymptotics for the Wiener sausage among Poissonian obstacles}
\author{Ryoki Fukushima
\footnote{Division of Mathematics, Graduate School of Science, Kyoto University, Kyoto, 606-8502, Japan; 
E-mail: fukusima@math.kyoto-u.ac.jp} }
\maketitle

\begin{abstract}
We consider the Wiener sausage among Poissonian obstacles. 
The obstacle is called $hard$ if Brownian motion entering the obstacle is 
immediately killed, and is called $soft$ if it is killed at certain rate. 
It is known that Brownian motion conditioned to survive among obstacles is confined in 
a ball near its starting point. 
We show the weak law of large numbers, large deviation principle in special cases and 
the moment asymptotics for the volume of the corresponding Wiener sausage. 
One of the consequence of our results is that the trajectory of Brownian motion 
almost fills the confinement ball. \\

\end{abstract} 

\section{Introduction}
We consider Brownian motion conditioned to avoid Poissonian obstacles. It is known that 
conditional Brownian motion typically localizes in a ball near its starting point under 
the annealed measure. 
In this article, we show that the trajectory of the particle almost 
fills the ball in which it is confined. 

We shall start by introducing the notation and the model. 
Let $\Omega$ be the set of locally finite simple pure point measures on $\R^d$ 
and $\p_{\nu}$ be the Poisson point process of constant intensity $\nu$ on $\Omega$. 
For a fixed nonpolar compact subset $K$ of $\R^d$ and $\Omega \ni \omega = 
\sum_i \delta_{x_i}$, we define the hard obstacles $S(\omega)=\bigcup_i (x_i+K)$.
Similarly, for a nonnegative, compactly supported and bounded measurable function $W$ 
which is not identically zero and 
$\Omega \ni \omega = \sum_i \delta_{x_i}$, we define the soft obstacles $V(x,\omega)
=\sum_i W(x-x_i)$. Next we introduce Brownian motion $Z_{\cdot}$ on $\R^d$ 
independent of the Poisson point process. The law of $Z_{\cdot}$ conditioned to start 
from $x \in \R^d$ is denoted by $P_{x}$ and $E_x$ stands for the corresponding expectation. 
For an open set $U \subset \R^d$ and a closed set $F \subset \R^d$, $T_U= \inf \set{s \ge 0 \,;
\, Z_s \notin U}$ and $H_F = \inf \set{s \ge 0 \,;\, Z_s \in F}$ are the exit time 
of $U$ and the entrance time of $F$, respectively.

We define the annealed path measure
\begin{equation*}
Q_t^{\mu, \nu}= \frac{1}{S_t^{\mu, \nu}}\exp \left\{-\int_0^t V(Z_s,\omega_1)ds \right\} 
1_{\{H_{S(\omega_2)}>t\}}\p_{\mu}^1 \otimes \p_{\nu}^2 \otimes P_0\label{1}
\end{equation*}
on $\Omega^2 \times C([0,t],\R^d)$ with $S_t^{\mu, \nu}$ the normalizing constant:
\begin{equation*}
S_t^{\mu, \nu}=\E_{\mu}^1\otimes \E_{\nu}^2\otimes E_0 \left[ \exp \left\{-\int_0^t 
V(Z_s,\omega_1)ds \right\} \, ; \, H_{S(\omega_2)}>t \right]\label{2}
\end{equation*}
with obvious notations.
This path measure describes the behavior of Brownian motion among the killing traps 
conditioned not to be killed up to time $t$. The first mathematical result concerned 
with this measure is Donsker-Varadhan's work \cite{DV75} about asymptotics for $S_t^{\nu,0}$ with $W$ as before 
and $S_t^{0, \nu}$ in the case $K$ is a closed ball of arbitrary fixed radius. They showed, 
using large deviation technique, that 
\begin{equation*}
S_t^{0, \nu}\;\left(\mathrm{or} \; S_t^{\nu,0}\right)=\exp \left\{ -c(d,\nu)t^{\frac{d}{d+2}}(1+o(1))\right\}\quad (t \to \infty)\label{3}
\end{equation*}
where
\begin{equation}
c(d,\nu) = \inf_{U \subset \R^d:\mathrm{open}} \set{\nu\left|U\right|+\lambda(U)}\\ \label{4}
\end{equation}
with $\left|U\right|$ the Lebesgue measure of $U$ and $\lambda(U)$ the principal Dirichlet 
eigenvalue of $-1/2 \Delta$ in $U$. It follows from Faber-Krahn's inequality
(see e.g.\ \cite{Ber03}) that balls with radius 
\begin{equation}
R_0(d,\nu)=\left( \frac{2\lambda_d}{d\nu\omega_d} \right)^{\frac{1}{d+2}}\label{5}
\end{equation}
achieve the infimum in \eqref{4}. Here $\omega_d$ is the volume 
of the $d$-dimensional unit ball $B(0,1)$ and $\lambda_d$  the principal Dirichlet 
eigenvalue of $-1/2 \Delta$ in $B(0,1)$. 
Therefore we can obtain the exact value of $c(d,\nu)$:
\begin{equation*}
 c(d,\nu)=\frac{d+2}{2}(\nu \omega_d)^{\frac{2}{d+2}}\Bigl(\frac{2\lambda_d}{d}\Bigr)^{\frac{d}{d+2}}.
\end{equation*}
Sznitman generalized this result to 
$S_t^{0, \nu}$ with arbitrary nonpolar compact $K$ and also improved the asymptotic estimates as follows:
\begin{equation}
 \begin{split}
  \exp \left\{ -c(d,\nu)t^{\frac{d}{d+2}}-\gamma(a,d,\nu)t^{\frac{d-1}{d+2}}\right\}
  &\le S_t^{0, \nu}\;\left(\mathrm{or} \; S_t^{\nu,0}\right)\\ 
  &\le \exp \left\{ -c(d,\nu)t^{\frac{d}{d+2}}+t^{\frac{d\mu(d)}{d+2}}\right\}\label{6}
 \end{split}
\end{equation}
for large $t$, using his `method of enlargement of obstacles'(see Theorem 4.5.6 in \cite{Szn98}). 
Here $\gamma(a,d,\nu)>0$ and $\mu(d) \in (0,1)$ are constants and $a$ is  
defined via
\begin{equation}
a=\inf \{u>0\,;\, K \;({\textrm{resp.\ }} \mathrm{supp}(W)) \subset \overline{B}(0,u)\}.\label{7}
\end{equation}
Sznitman($d=2$, in \cite{Szn91b}) and Povel($d \ge 3$, in \cite{Pov99}), motivated 
by the proof of the lower bound in \eqref{6}, showed that surviving Brownian particle is 
typically confined in a ball with radius $t^{1/(d+2)}R_0$ for large $t$. 
\begin{thm}\label{Theorem1}(Confinement property)\\
Let $d \ge 2$. There exist constants $\kappa_1 >1$ and $0<\kappa_2<1$ and for each 
$(\omega_1,\omega_2) \in \Omega^2$ a ball $B(\omega_1,\omega_2)$ with center in 
$B(0, R_0(d,\mu+\nu)+\kappa_1 t^{-\kappa_2/(d+2)})$ and with radius in 
$[R_0(d,\mu+\nu),R_0(d,\mu+\nu)+\kappa_1 t^{-\kappa_2/(d+2)}]$ 
such that
\begin{equation*}
\lim_{t \to \infty}Q_t^{\mu,\nu}\left(T_{t^{1/(d+2)}B(\omega_1,\omega_2)}>t\right)=1.\label{8}
\end{equation*}
\end{thm}\ \\

\noindent
Although Sznitman and Povel showed this theorem only in the case $Q_t^{0, \nu}$, 
their argument is easily applicable to above version. 
As a consequence of this property, the volume of 
the Wiener sausage $W_t^C= \bigcup_{s \le t}(Z_s+C)$ associate with a compact set 
$C \subset \R^d$ is typically not larger than 
$t^{d/(d+2)}\left|B(0,R_0)\right|$ under $Q_t^{\mu, \nu}$. 

The first result of this paper is that $|W_t^C|$ under $Q_t^{\mu, \nu}$ asymptotically 
equals to $t^{d/(d+2)}\left|B(0,R_0)\right|$ in the sense of the weak 
law of large numbers:
\begin{thm}\label{Theorem2}
Let $d \ge 2$, $\mu \ge0$, $\nu \ge 0$ and $\mu + \nu >0$. Then we have for any 
nonpolar compact set $C$ and $\epsilon>0$,
\begin{equation}
\lim_{t \to \infty}Q_t^{\mu,\nu}\left(\left|t^{-\frac{d}{d+2}}|W_t^C|-\omega_d 
R_0(d,\mu+\nu)^d \right|>\epsilon \right)=0.\label{9}
\end{equation}
Moreover, if $\mu=0$, $\nu>0$ and $C \subset -K$, the law of $t^{-d/(d+2)}|W_t^C|$ under
$Q_t^{0,\nu}$ satisfies following large deviation principle:
\begin{equation}
 \begin{split}
-\inf_{x \in \Gamma^{\circ}} I(x) \le
&\liminf_{t \to \infty} t^{-\frac{d}{d+2}} \log Q_t^{0,\nu}\left(t^{-\frac{d}{d+2}}|W_t^C| 
\in \Gamma \right) \\ \le
&\limsup_{t \to \infty} t^{-\frac{d}{d+2}} \log Q_t^{0,\nu}\left(t^{-\frac{d}{d+2}}|W_t^C| 
\in \Gamma \right) \le 
-\inf_{x \in \overline{\Gamma}} I(x),\label{10}
 \end{split}
\end{equation}
where $\Gamma$ is arbitrary Borel subset of $(0,\infty)$ and rate function $I$ is given by 
\begin{equation*}
I(x)=\nu x + \lambda_d \left(\frac{\omega_d}{x}\right)^{\frac{2}{d}}-c(d,\nu).\label{11}
\end{equation*}
\end{thm}\ \\

\noindent{\bf Remark.} The assumption $C \subset -K$ may look rather technical. 
But the large deviation principle with above rate function fails when $C$ is much larger than $K$. 
We shall give an example after the proof of Theorem 2. \\

\noindent
Theorem \ref{Theorem2}, combined with Theorem \ref{Theorem1}, implies that the Wiener sausage under 
$Q_t^{\mu,\nu}$ covers almost all area of the ball in which it is confined. 

The next result is the improvement of the convergence to $L^p$ sense. 
We derive it as a corollary of following exponential tightness estimate:
\begin{thm}\label{Theorem3}
Let $d \ge 2$, $\mu \ge0$, $\nu \ge 0$, $\mu + \nu >0$. Then for any $\eta >0$, 
\begin{equation*}
\sup_{t\ge 1} Q_t^{\mu,\nu}\left(\exp\left\{\eta t^{-\frac{1}{d+2}}\sup_{0\le s \le t}
\left|Z_t \right|\right\}\right)<\infty.\label{12}
\end{equation*}
\end{thm}\ \\

\begin{cor}\label{Corollary}
Under the same conditions as in Theorem \ref{Theorem3}, we have 
\begin{equation}
\sup_{t\ge 1} Q_t^{\mu,\nu}\left(\exp\left\{\eta\left(t^{-\frac{d}{d+2}}|W_t^C|\right)
^{\frac{1}{d}}\right\}\right)<\infty \label{12'}
\end{equation}
for any $\eta>0$. Consequently, we have for all $p>0$
\begin{equation}
\lim_{t \to \infty}Q_t^{\mu,\nu}\left(\left|t^{-\frac{d}{d+2}}|W_t^C|
-\omega_d R_0(d,\mu+\nu)^d \right|^p\right)=0.\label{12''}
\end{equation}
\end{cor}
\noindent{\it Proof.} 
Since $W_t^C \subset B(0,\sup_{0\le s\le t}|Z_s|+\mathrm{diam}C)$, \eqref{12'} follows directly from \eqref{12}. 
From \eqref{12'}, we have that $\{t^{-dp/(d+2)}|W_t^C|^p\}_{t\ge 0}$ is 
uniformly integrable for any $p>0$, which implies \eqref{12''}.
\hfill $\square$\\

Now let us briefly explain the construction of this article. We start by considering 
exponential moments of $|W_t^C|$. Since negative exponential moments are easy 
to estimate, lower estimate of \eqref{9} follows from rather simple calculations. 
In the special case $\mu=0$, $\nu>0$ and $C \subset -K$, the upper bound for 
exponential moments can be extended to positive parameters. Then, we can derive large deviation 
upper bound using a similar argument to the G\"{a}rtner-Ellis theorem. 
The large deviation lower bound is obtained by considering a specific strategy for 
Wiener sausage to achieve given volume. 
Next, we shall give the proof of Theorem 3. 
Our strategy is essentially the same as the Povel's proof of Theorem \ref{Theorem1} 
but we need quantitatively refined estimate for the probability of the process exiting the confinement ball. 
The proofs of Theorem \ref{Theorem1} and the upper estimate of \eqref{9} will also be given 
along the way in order to make this article reasonably self-contained. 

\section{Lower estimate of Theorem \ref{Theorem2} and large deviation}
In this section, we are going to show the lower estimate of Theorem \ref{Theorem2} and 
the large deviation result. Firstly, note that we can prove 
\begin{equation}
S_t^{\mu, \nu} 
=\exp \left\{ -c(d,\mu+\nu)t^{\frac{d}{d+2}}(1+o(1))\right\}\label{14}   
\end{equation}
when $t$ goes to $\infty$. 
Indeed, the lower bound is obvious since $S_t^{\mu,\nu}$ is bounded from below by $S_t^{0,\mu+\nu}$ 
with hard obstacle $({\rm supp}\,W) \cup K$. And the upper bound follows from 
the same argument as in the proof of theorem 4.5.6 in \cite{Szn98}, 
using the method of enlargement of obstacles which will be explained in section 3.1. 
(In the upper bound of theorem 4.5.6 in \cite{Szn98}, the `enlarged obstacles' are 
mainly considered and therefore the shape of obstacles has little to do with the argument.) 

Similarly, for all $\lambda<0$ and large enough $t$ we can prove 
\begin{equation}
 \begin{split} 
&\E_{\mu}^1\otimes \E_{\nu}^2\otimes E_0 \left[ \exp \left\{-\int_0^t V(Z_s,\omega_1)ds
  +\lambda|W_t^C |\right\}\, ; \, H_{S(\omega_2)}>t \right]\\
=\,&\E_{\mu}^1\otimes \E_{\nu}^2\otimes \E_{-\lambda}^3\otimes E_0 \left[ \exp \left\{-\int_0^t V(Z_s,\omega_1)ds
  \right\}\, ; \, H_{S(\omega_2)\cup \tilde{S}(\omega_3)}>t \right]\\
=\,&\exp \left\{ -c(d,\mu+\nu-\lambda)t^{\frac{d}{d+2}}(1+o(1))\right\},\label{14'}
 \end{split}
\end{equation}
where we set $\tilde{S}(\omega)=\bigcup_i(x_i-C)$ for $\omega=\sum_i \delta_{x_i}$.
As a consequence, 
\begin{equation}
\begin{split}
&Q_t^{\mu,\nu}\left(\exp\left\{\lambda |W_t^C |\right\}  \right)\\
=\, &\exp \left\{ \left(c(d,\mu+\nu)-c(d,\mu+\nu-\lambda)\right)
t^{\frac{d}{d+2}}(1+o(1))\right\}\label{15}
\end{split}
\end{equation}
for $\lambda<0$ when $t \to \infty$. 
Here we have implicitly used the fact $c(d,\mu+\nu) \neq c(d,\mu+\nu-\lambda)$ to ensure that the $o(1)$ 
in \eqref{14} and \eqref{14'} is again $o(1)$ in \eqref{15}. 
Now we prove the lower estimate of \eqref{9} in Theorem 2. 

\begin{prop} For any $\epsilon > 0$, 
\begin{equation}
\lim_{t \to \infty}Q_t^{\mu,\nu}\left(t^{-\frac{d}{d+2}}|W_t^C| < \omega_d R_0(d,\mu+\nu)^d - \epsilon \right)=0.
\end{equation}
\end{prop}
{\it{Proof.}} 
Let $m = \omega_d R_0(d,\mu+\nu)^d - \epsilon$. 
Then Chebyshev's inequality and \eqref{15} yield 
\begin{equation}
 \begin{split}
&Q_t^{\mu,\nu}\left(t^{-\frac{d}{d+2}}|W_t^C| \le m \right)\\
\le &\exp \left\{ \left(-\lambda m+c(d,\mu+\nu)-c(d,\mu+\nu-\lambda)\right)
t^{\frac{d}{d+2}}(1+o(1))\right\} \\
\le &\exp \left\{ \lambda\left(\frac{c(d,\mu+\nu)-c(d,\mu+\nu-\lambda)}{\lambda}-m\right)
t^{\frac{d}{d+2}}(1+o(1))\right\}\label{16}
 \end{split}
\end{equation}
as $t \to \infty$, provided $(c(d,\mu+\nu)-c(d,\mu+\nu-\lambda))/\lambda-m \neq 0$. 
If we note the fact that 
\begin{equation}
\omega_d R_0(d,\nu)^d
=\omega_d^{\frac{2}{d+2}}\nu^{-\frac{d}{d+2}}\Bigl(\frac{2\lambda_d}{d}\Bigr)^{\frac{d}{d+2}}
=\frac{\partial}{\partial \nu}c(d,\nu), \label{17}
\end{equation}
we can actually find a $\lambda<0$ such that $(c(d,\mu+\nu)-c(d,\mu+\nu-\lambda))/\lambda-m > 0$. 
Since the right hand side of \eqref{16} goes to 0 for this $\lambda$, we have done. 
\hfill $\square$\\

Next, we shall prove large deviation result. 
We start with the upper bound. 
\begin{prop} Suppose $\mu=0$, $\nu>0$, and $C \subset -K$. 
Then for arbitrary Borel subset $\Gamma \subset (0,\infty)$, 
\begin{equation*}
\limsup_{t \to \infty} t^{-\frac{d}{d+2}} \log Q_t^{0,\nu}\left(t^{-\frac{d}{d+2}}|W_t^C| 
\in \Gamma \right) \le -\inf_{x \in \overline{\Gamma}} I(x). 
\end{equation*}
Here $I$ is the rate function defined in Theorem 2. 
\end{prop}
{\it{Proof.}} 
By the assumption $C \subset -K$, we can extend the upper bound of \eqref{15} 
to $0 < \lambda \le \nu$ as follows:
\begin{equation*}
 \begin{split}
Q_t^{0,\nu}\left(\exp\left\{\lambda |W_t^C |\right\} \right)
\le &\, \frac{E_0 \left[ \exp \left\{\lambda|W_t^{-K} | -\nu|W_t^{-K} | 
 \right\}\right]}{E_0 \left[ \exp \left\{-\nu|W_t^{-K} |  \right\}\right]}\\
=&\, \exp \left\{ \left(c(d,\nu)-c(d,\nu-\lambda)\right) t^{\frac{d}{d+2}}(1+o(1)) \right\}.\label{20}
 \end{split}
\end{equation*}
Here we have used $c(d,\nu)-c(d,\nu-\lambda) \neq 0$ as in the derivation of \eqref{15}. 
Therefore, we have following upper bound on the logarithmic generateing function: 
\begin{equation}
\begin{split}
&\limsup_{t \to \infty} t^{-\frac{d}{d+2}} \log Q_t^{0,\nu}
 \left(\exp\left\{\lambda |W_t^C |\right\} \right)\\
\le &\,\left\{
\begin{array}{lr}
 c(d,\nu)-c(d,\nu-\lambda) &(\lambda \le \nu),\\[8pt]
 \infty  &(\lambda > \nu).
\end{array}\right.\label{22}
\end{split}
\end{equation} 
Then, the large deviation 
upper bound follows from very similar argument to the proof of the G\"{a}rtner-Ellis 
theorem (cf. \cite{DZ98}) and 
the rate function is given by the Fenchel-Legendre transform of the right hand side of \eqref{22}. 
\hfill $\square$\\

Next, we go on to the lower bound, which do not require the assumption $C \subset -K$. 
\begin{prop} Suppose $\mu=0$ and $\nu>0$. 
Then for arbitrary Borel subset $\Gamma \subset (0,\infty)$, 
\begin{equation*}
\liminf_{t \to \infty} t^{-\frac{d}{d+2}} \log Q_t^{0,\nu}\left(t^{-\frac{d}{d+2}}|W_t^C| 
\in \Gamma \right) \ge -\inf_{x \in \Gamma^{\circ}} I(x).
\end{equation*}
Here $I$ is the rate function defined in Theorem 2. 
\end{prop}
{\it{Proof.}} 
It is enough to prove 
\begin{equation}
\liminf_{t \to \infty} t^{-\frac{d}{d+2}} \log Q_t^{0,\nu}\left(t^{-\frac{d}{d+2}}|W_t^C| 
\in (x-\delta,x+\delta) \right) \ge -I(x)\label{23'}
\end{equation}
for any $x>0$ and $\delta>0$. To this end, set $r(t)=t^{1/(d+2)}(x/\omega_d)^{1/d}$ and consider 
the specific event 
\begin{equation*}
A_1 \times A_2=\set{\omega\left(B(0,r(t)+a)\right)=0} \times \set{T_{B(0,r(t))}>t}
\end{equation*}
where $a$ was defined in \eqref{7}. (Note that $C \subset -K \subset \overline{B}(0,a)$.) Since we know 
\begin{equation*}
 \begin{split}
&Q_t^{0,\nu}(A_1 \times A_2) = \frac{1}{S_t^{0,\nu}}\p_{\nu}(A_1) P_0(A_2)\\
\ge\, & \mathrm{const}(d)\exp\left\{-\left(\nu \omega_d r(t)^d+\frac{\lambda_d}{r(t)^2}t-c(d,\nu)
t^{\frac{d}{d+2}}\right)(1+o(1))\right\}\\
=\, & \mathrm{const}(d)\exp\left\{-I(x)t^{\frac{d}{d+2}}(1+o(1))\right\}
 \end{split}
\end{equation*}
from \eqref{6} and a well known eigenfunction expansion, \eqref{23'} follows once we have shown
\begin{equation}
\lim_{t \to \infty}P_0\left(t^{-\frac{d}{d+2}}|W_t^C| \in (x-\delta,x+\delta)\,\Bigl|\, A_2 \right)=1.\label{23''}
\end{equation}
\textit{Proof of \eqref{23''}.} Using Brownian scaling by the scale $\epsilon=t^{-1/(d+2)}$, \eqref{23''} is equivalent to
\begin{equation}
\lim_{t \to \infty}P_0\left(|W_{t\epsilon^2}^{\epsilon C}| \in (x-\delta,x+\delta)\,
\Bigl|\, T_{B(0,r(1))}>t\epsilon^2 \right)=1.\label{23-1}
\end{equation}
To show this we shall use Theorem 3.2.3 in \cite{Szn98}, which claims 
\begin{equation*}
\lambda(U \setminus K)\ge \lambda(U)+\frac{\mu(U)-\lambda(U)}{\mu(U)}\inf_K \varphi_U^2 \cdot\mathrm{cap}_U(K)
\end{equation*}
for any nonempty bounded domain $U \subset \R^d$ and compact set $K \subset U$. Here $\varphi_U$ denotes the 
positive $L^2$-normalized principal 
eigenfunction and $\mu(U)$ the second smallest eigenvalue both associated with $-\Delta /2$ in $U$ with Dirichlet 
boundary condition. From this and the fact that $\varphi_{B(0,r)}$ has nondegenerate gradient near the boundary, 
we have for all $x \in B(0,r(1)-\sqrt{\epsilon})$, 
\begin{equation}
 \begin{split}
&\log P_0\left(H_{x-\epsilon C}>t\epsilon^2 \,\Bigl|\, T_{B(0,r(1))}>t\epsilon^2 \right)\\
= &-\left(\lambda(B(0,r(1))\setminus (x-\epsilon C))-\lambda(B(0,r(1))) \right)t^{\frac{d}{d+2}}(1+o(1))\\
\le &-\mathrm{const}(d,r(1))\epsilon\cdot \mathrm{cap}_{B(0,r(1))}(\epsilon C)t^{\frac{d}{d+2}}(1+o(1))\label{23'''}
 \end{split}
\end{equation}
as $t \to \infty$. The right hand side of \eqref{23'''} goes to $-\infty$ as $t \to \infty$ since 
\begin{equation*}
\mathrm{cap}_{B(0,r(1))}(\epsilon C) \sim \left\{
\begin{array}{lr}
 \mathrm{const}(C,2)\left(\log\frac{1}{\epsilon}\right)^{-1} &(d=2),\\[8pt]
 \mathrm{const}(C,d)\epsilon^{d-2}  &(d \ge 3).
\end{array}\right.
\end{equation*}
As a consequence, we have
\begin{equation*}
 P_0\left(H_{x-\epsilon C}>t\epsilon^2 \,\Bigl|\, T_{B(0,r(1))}>t\epsilon^2 \right) \to 1 \quad (t \to \infty)
\end{equation*}
and therefore
\begin{equation*}
 \begin{split}
&E_0\left[|W_{t\epsilon^2}^{\epsilon C}| \,\Bigl|\, T_{B(0,r(1))}>t\epsilon^2\right]\\
=\,&\int_{\R^d}P_0\left(H_{x-\epsilon C}\le t\epsilon^2 \,\Bigl|\, T_{B(0,r(1))}>t\epsilon^2 \right)dx\\
\ge \,&\int_{B(0,r(1)-\sqrt{\epsilon})}P_0\left(H_{x-\epsilon C}\le t\epsilon^2 \,\Bigl|\, T_{B(0,r(1))}>t\epsilon^2 \right)dx\\
\to \,&|B(0,r(1))|=x \qquad (t \to \infty).
 \end{split}
\end{equation*}
This, together with the obvious fact that 
\begin{equation*}
P_0\left(|W_{t\epsilon^2}^{\epsilon C}| \le |B(0,r(1)+a\epsilon)| \,\Bigl|\, T_{B(0,r(1))}>t\epsilon^2 \right)=1,
\end{equation*}
implies \eqref{23-1}.
\hfill $\square$\\
Now that we have shown \eqref{23''}, the proof of Proposition 3 is completed. \hfill $\square$\\

Finally, we shall give an example noticed after Theorem 2. 
Basically, it comes from the case where the `ballistic strategy' dominates above `localizing strategy'. 
{\flushleft{\bf Example.}} 
Let $C=\overline{B}(0,R)$ and $K=\overline{B}(0,1)$ where $R$ will be taken large. 
The key is to consider the specific strategy:
\begin{equation*}
 A= \set{Z_ t \in B\left(t^{\frac{d}{d+2}}{h}, rt^{\frac{d}{d+2}}\right) } 
\end{equation*}
where $h \in \R^d$ and $r>0$ satisfy $|h|>r$. 
On this event, we have
\begin{equation}
 |W_t^C|>t^{d/(d+2)}(|h|-r)\omega_{d-1}R^{d-1} \label{ex1}
\end{equation}
by considering cross sections orthogonal to $h$. 
For the `cost' of this strategy, we use the large deviation estimate
\begin{equation}
 t^{-\frac{d}{d+2}}\log Q_t^{0,\nu}(A) \sim -\inf_{x \in B(h,r)}\beta_0(x)
 \quad (t \to \infty)\label{ex2}
\end{equation}
which was shown by Sznitman in \cite{Szn95a}. Here $\beta_0$ is the annealed Lyapunov exponent introduced in \cite{Szn95a}, 
which measures the decay rate of the probability for Brownian motion to perform a long crossing among $S(\omega)$. 

Combining \eqref{ex1} and \eqref{ex2}, we have 
\begin{equation*}
 \begin{split}
  &\liminf_{t \to \infty}t^{-\frac{d}{d+2}}\log Q_t^{0,\nu}\left(t^{-d/(d+2)}|W_t^C|>(|h|-r)\omega_{d-1}R^{d-1}\right) \\
  &\ge \,-\inf_{x \in B(h,r)}\beta_0(x). \label{ex3}
 \end{split}
\end{equation*} 
Since the right hand side is independent of $R$, the upper bound of \eqref{10} breaks down when $R$ is large. 

\section{Upper estimates}
We shall prove the upper estimate of \eqref{9} and Theorem \ref{Theorem3} in this section. 
As described at the end of the section 1, the proofs of them are 
based on Theorem \ref{Theorem1} and its proof. The starting point is to adopt the scale
\begin{equation*}
\epsilon = t^{-\frac{1}{d+2}}\label{24}
\end{equation*}
and consider $\epsilon^{-1} Z_{t\epsilon^2 }$, $\p^1_{\mu\epsilon^{-d}}$ and 
$\p^2_{\nu\epsilon^{-d}}$. 
We introduce the notation ${\bold E}_{\epsilon}=\E^1_{\mu \epsilon^{-d}} \otimes \E^2_{\nu \epsilon^{-d}}$ for simplicity. 
Then, for instance, Theorem \ref{Theorem1} follows once 
we have shown that for all $(\omega_1,\omega_2) \in \Omega^2$ there exists a ball 
$B(\omega_1,\omega_2)$ with radius in $[R_0,R_0+\kappa_1 \epsilon^{\kappa_2}]$ and 
center in $B(0,R_0+\kappa_1 \epsilon^{\kappa_2})$ such that 
\begin{equation*}
 \begin{split}
\lim_{t \to \infty} \frac{1}{S_t^{\mu,\nu}} {\bold E}_{\epsilon} \otimes E_0 \biggl[ \exp &\left\{-\int_{0}^{\tau}
V_{\epsilon}(Z_s, \omega_1)ds \right\}\,; \\
&T_{B(\omega_1,\omega_2)} \wedge H_{S_{\epsilon}(\omega_2)}>\tau \biggr]=1.\label{25}
 \end{split}
\end{equation*}
Here $\tau= t\epsilon^2$ and
\begin{gather*}
V_{\epsilon}(x, \omega)=\epsilon^{-2} \sum_i W\left(\frac{x-x_i}{\epsilon}\right),\label{26}\\
S_{\epsilon}(\omega)=\bigcup_i (x_i+\epsilon K)\label{27}
\end{gather*}
for $\omega=\sum_i \delta_{x_i}$. For $(\omega_1,\omega_2) \in \Omega^2$ and 
open set $U$ we also define
\begin{equation*}
\lambda^{\epsilon}_{\omega_1,\omega_2}(U)=\lambda_{V_{\epsilon}(\cdot,\omega_1)}(U\setminus
S_{\epsilon}(\omega_2))\label{28}
\end{equation*}
where $\lambda_V(U)$ denotes the principal Dirichlet eigenvalue of $-1/2\Delta+V$ in $U$.

\subsection{Method of enlargement of obstacles}
In this subsection, we shall recall elements and some estimates from the method 
of enlargement of obstacles in \cite{Szn98}. 
The method is based on coarse graining of the space and 
construction of two disjoint sets $\D$ and $\B$ 
for $\omega \in \Omega$ and $\epsilon \in (0,1)$. The set $\D$ 
is called `density set', where one enlarge (both the support and the height of) 
the obstacles and the set $\B$ is called `bad set', where the 
obstacles exist but are left almost untouched. 

To construct these sets, we need parameters $0 < \alpha < \gamma < \beta <1$, 
$\delta >0$ and an integer $L \ge 2$. Using these parameters, we introduce 
spatial scales $1 \gg \epsilon^{\alpha} \gg \epsilon^{\gamma} \gg \epsilon^{\beta} 
\gg \epsilon$ and $L$-adic decomposition of $\R^d$. 
We also need following notation concerning $L$-adic decomposition of $\R^d$. 
Let $\mathcal{I}_k$ be the collection of indices of the form
\begin{equation*}
 \index= (i_0, i_1, \ldots , i_k) \in \Z^d \times (\set{0, 1,\ldots,L-1}^d)^k. 
\end{equation*}
We associate to above index $\index$ a box: 
\begin{equation*}
 C_{\index} = q_{\index} + L^{-k}[0,1]^d \;\textrm{ where }\; q_{\index} = i_0 + L^{-1}i_1 + \cdots +L^{-k}i_k. 
\end{equation*}
For $\index\in \mathcal{I}_k$ and $k' \le k$, we define the truncation 
\begin{equation*}
 [\index]_{k'} = (i_0, i_1, \ldots , i_{k'}). 
\end{equation*}
Finally, we pick integers $n_{\alpha}
(\epsilon)$, $n_{\gamma}(\epsilon)$ and $n_{\beta}(\epsilon)$ such that
\begin{equation}
L^{-n_{*}-1} \le \epsilon^{*} < L^{-n_{*}}\quad (\textrm{$*$ is $\alpha$, 
$\gamma$ or $\beta$}).\label{29}
\end{equation}
$L^{-n_{*}}$ plays the role of the scale $\epsilon^{*}$ in the context of $L$-adic 
decomposition. 

For $\omega = \sum_q \delta_{x_q}$ and $\index\in \mathcal{I}_k$, 
\begin{equation}
 K_{\index} = L^k \bigg( \bigcup_{x_q \in C_{\index}} \overline{B}(x_q, a\epsilon) \bigg) \label{skelton}
\end{equation} 
is called the skeleton of traps. Here $a$ is taken so large as $\overline{B}(0, a)$ includes both 
$K$ and ${\rm supp}\,W$. 
(We take larger $a$ so that $\overline{B}(0, a)$ includes $-C$ for the proof of \eqref{14'}.) 
Using this skeleton, the density set $\D$ is defined as follows:
\begin{Def} ((4.2.13) in \cite{Szn98}) 
 $C_{\index}$ $(\index\in \mathcal{I}_{n_{\gamma}})$ is called a density box if it satisfies the 
 quantitative Wiener criterion:  
 \begin{equation*}
  \sum_{n_{\alpha} < k \le n_{\gamma}} {\rm cap}(K_{[\index]_k}) \ge \delta (n_{\gamma}-n_{\alpha}). \label{qWc}
 \end{equation*}
 Here ${\rm cap}(\,\cdot\,)$ denotes the capacity relative to $1-\Delta/2$ when $d=2$ 
 and $-\Delta/2$ when $d \ge 3$. 
 The union of all density boxes is denoted by ${\D}$. 
\end{Def}
Next, the bad set $\B$ is defined as follows: 
\begin{Def} ((4.3.47) in \cite{Szn98}) 
 $C_{\index}$ $(\index\in \mathcal{I}_{n_{\beta}})$ is called a bad box if 
 $\omega (C_{\index}) \ge 1$ and $C_{\index} \not \subset \D$. 
 The union of all bad boxes is denoted by $\B$. 
\end{Def}

As a result of above construction, $\D$ and $\B$ satisfy
\begin{gather*}
\D \cap \B = \emptyset,\label{30}\\
\omega \left(\R^d \setminus \left(\D \cup \B \right) \right)=0,\label{31}
\end{gather*}
\begin{equation*}
 \begin{split}
&\textrm{for each box }C_q=q+[0,1)^d,\; q \in \Z^d, \textrm{ the sets }\D \cap C_q\\
&(\textrm{resp. } \B \cap C_q)\textrm{ can take no more than }2^{\epsilon^{-d\gamma}}
(\textrm{resp. }2^{\epsilon^{-d\beta}})\\
&\textrm{different shapes as } \omega \textrm{ varies over } \Omega.\label{32}
 \end{split}
\end{equation*}
The notation and definitions are exactly the same as in \cite{Szn98} so far. 

Next, we state four estimates which is the mixed obstacles version of the results in \cite{Szn98}. 
For the proofs of these estimates, we shall put some comments at the end of this subsection. 
We define the density set $\den$ and the bad set $\bad$ as above by letting $\omega=\omega_1+\omega_2$. 
The first claims that solidifying $\clden$, i.e.\ imposing Dirichlet conditions on $\clden$,
does not cause essential increase of the principal eigenvalues.\\

\noindent\textbf{Spectral control I} (Theorem 4.2.3 in \cite{Szn98}) 
There exist $c_1(d,W,K)>0$ such that for all 
$\rho \in (0, \delta c_1 \frac{\gamma-\alpha}{(d+2)\log L})$ and $M>0$, 
\begin{equation*}
\lim_{\epsilon \to 0} \sup_{(\omega_1,\omega_2) \in \Omega^2, \,U:\mathrm{open}}\epsilon^{-\rho}
\left(\lambda^{\epsilon}_{\omega_1,\omega_2}\left(U \setminus \clden \right)\wedge M 
-\lambda^{\epsilon}_{\omega_1,\omega_2}\left(U \right)\wedge M \right)=0.\label{33}
\end{equation*}\\

The second corresponds to the volume of the bad set. Since we cannot 
control the solidifying effect on $\bad$, we need to show that it is not too large. \\

\noindent\textbf{Volume control} (Theorem 4.3.6 in \cite{Szn98}) 
There exist $L \ge 2$, $\delta>0$ and $\kappa>0$ such that
\begin{equation*}
\lim_{\epsilon \to 0} \sup_{\omega \in \Omega,\,q \in \Z^d}\epsilon^{-\kappa}
|\bad \cap C_q| < 1.\label{34}
\end{equation*}\\

The third estimate says that the region where $\den^c$ is locally thin is hard to survive. 
Therefore one can expect that, as in the first estimate, solidifying such a region does not 
cause essential increase of the principal eigenvalues. This is precisely the role 
of the fourth estimate. \\

\noindent\textbf{Spectral control II} (Proposition 4.2.4 in \cite{Szn98}) 
There exist $c_2(d)>0$ such that for all $\epsilon \in (0,1)$, $r \in (0,1/4)$, 
$(\omega_1, \omega_2) \in \Omega^2$ and open set $U$ satisfying:
\begin{gather*}
4a\epsilon<L^{-n_{\gamma}(\epsilon)}<L^{-n_{\alpha}(\epsilon)}<r,\label{35}\\
\delta c_1(n_{\gamma}(\epsilon)-n_{\alpha}(\epsilon))>\log 2,\label{36}\\
\sup_{q \in \Z^d} \left| \left(U\setminus \clden \right)\cap C_q \right| < r^d,\label{37}
\end{gather*}
one has 
\begin{equation*}
\lambda^{\epsilon}_{\omega_1,\omega_2}(U)>\frac{c_2}{r^2}.\label{38}
\end{equation*}\\

\noindent\textbf{Spectral control III} (Theorem 4.2.6 in \cite{Szn98}) 
For all $M>0$, there exist $c_3(d)>0$, $c_4(d,M)>1$ and $r_0(d,M) \in (0,1/4)$ 
such that
\begin{equation}
\limsup_{\epsilon \to 0}\, \widetilde{\sup}\,
\exp \left\{ c_3\left[ \frac{R}{4r} \right]\right\} 
\left(\lambda^{\epsilon}_{\omega_1,\omega_2}\left(U_1 \right)\wedge M 
-\lambda^{\epsilon}_{\omega_1,\omega_2}\left(U_2 \right)\wedge M \right) \le 1,\label{39}
\end{equation}
where $[\,\cdot\,]$ denotes the integer part and $\widetilde{\sup}$ the supremum 
over all $(\omega_1, \omega_2) \in \Omega^2$, $U_1 \subset U_2\,:\,$open, $R>0$ and $r>0$ 
such that for some closed set $A$,
\begin{gather}
L^{-n_{\alpha}(\epsilon)}<r<r_0,\label{40}\notag \\
\frac{R}{4r}>c_4,\label{41}\notag \\
\sup_{q \in \Z^d} \left| \left(U_2\setminus \left( A \cup \clden \right) \right)
\cap C_q \right| < r^d,\label{42}\\
\mathrm{dist}_{\|\cdot \|}\left(U_2 \setminus U_1, A \cap U_2 \right) \ge R.\label{43}
\end{gather}
Here, $\| \cdot \|$ is the maximal norm on $\R^d$.\\

For a typical case \eqref{39} applies, we introduce 
\begin{gather*}
\A_{\epsilon}(\omega_1,\omega_2)=\bigcup_{C_q\,:\,|C_q \setminus \den| \ge r^d}\overline{C}_q,\label{44}\\
\Ocl=\set{x \in \R^d;\,\mathrm{dist}_{\|\cdot \|}(x,\A_{\epsilon}(\omega_1,\omega_2))<R}.\label{45}
\end{gather*}
Each $C_q$ contained in $\A_{\epsilon}(\omega_1,\omega_2)$ is called \textit{clearing box} 
and $\A_{\epsilon}(\omega_1,\omega_2)$ is called \textit{clearing set}. 
Then for any open set $\T$ and $R>0$, 
\[
U_1=\T \cap \Ocl, \quad U_2=\T
\]
satisfy \eqref{42} and \eqref{43} with $A=\A_{\epsilon}(\omega_1,\omega_2)$. 

Before closing this subsection, we briefly explain how to prove these estimates. 
The volume estmate is equivalent to that in \cite{Szn98} since our bad set is the same as 
that for hard obstacles with $K=\overline{B}(0,a)$. 
For the spectral controls, we mention that in the proofs in \cite{Szn98}, 
the dependence on the shapes of obstacles only appears in following key lemma: \\

\noindent\textbf{Key Lemma} (Lemma 4.2.1 in \cite{Szn98})
 There exist $c_1(d,W)>0$ (or $c_1(d,K)$ for hard obstacles) such that if $4a\epsilon < L^{-n_{\gamma}}$,  
 \begin{gather*}
  E_x\left[ \exp\set{-\int_0^{H_{n_{\alpha}}(\epsilon)} V_{\epsilon}(Z_s, \omega) ds } \right]
  \le \exp\set{-c_1 \sum_{n_{\alpha} < k \le n_{\gamma}} {\rm cap}(K_{[\i]_k})}\\
  \left({\rm resp.\ } P_x(H_{n_{\alpha}} < H_{S(\omega)})
  \le \exp\set{-c_1 \sum_{n_{\alpha} < k \le n_{\gamma}} {\rm cap}(K_{[\i]_k})}\right)
 \end{gather*}
 for any $\omega \in \Omega$, $\index\in \mathcal{I}_{n_{\gamma}}$ and $x \in C_{\index}$. 
 Here $H_{n_{\alpha}}=\inf\{s \ge 0\,;\, \|Z_s -Z_0 \| \ge L^{-n_{\alpha}}\}$.\\

\noindent
However, it is routine to extend this key lemma to our mixed obstacles 
if we replace $c_1$ by $c_1(d,W) \wedge c_1(d,K)$. (This is the constant $c_1(d,W,K)$ appeared in spectral control I and II.) 
Thus we can prove the spectral controls in our setting by exactly the same ways as in \cite{Szn98}.  
 
\subsection{Construction of the confinement ball}
In this subsection, we shall construct $B(\omega_1,\omega_2)$ in Theorem \ref{Theorem1}. 
Since results are essentially the same as in \cite{Pov99}, we omit the proofs and refer 
counterparts instead. 
From now on, we fix an admissible collection of parameters 
\begin{equation*}
\alpha,\; \beta,\; \gamma,\; \delta,\; L,\; \rho,\; \kappa\label{46} 
\end{equation*}
and pick $R=1$ and
\begin{equation*}
r=\epsilon^{\alpha_0}, \textrm{ with } 0< \alpha_0 < \min \left(\alpha, 1-\beta, 
\frac{\kappa}{d}\right)<1\label{47}
\end{equation*}
which allow us to apply the results in section 3.1.
For an explanation about the admissible collection of parameters, we refer reader to 
the remarks after (4.3.66) in \cite{Szn98}. 

Now let us start by introducing the open set
\begin{equation*}
\T=(-[t],[t])^d\label{48}
\end{equation*}
and notation
\begin{equation*}
v(\tau)=\exp \left\{-\int_{0}^{\tau}V_{\epsilon}(Z_s, \omega_1)ds \right\}
1_{\{H_{S_{\epsilon}(\omega_2)}>\tau\}}.\label{49}
\end{equation*}
Then, using standard estimates on Brownian motion and \eqref{14}, we have 
\begin{equation}
 \begin{split}
\frac{1}{t}\log\frac{1}{S_t^{\mu,\nu}}{\bold E}_{\epsilon} \otimes E_0 \left[ v(\tau) \,;\,T_{\T}\le \tau\right]<0\label{50}
 \end{split}
\end{equation}
for large $t$. Since this is good enough for our purpose, we restrict our consideration on 
$\{T_{\T}>\tau\}$ in the sequel. We also introduce the open set 
\begin{equation*}
\mathscr{U}(\omega_1,\omega_2)=(\T \cap \Oc)\setminus \clden.\label{51}
\end{equation*}
Then, we have following constraint on this set.
\begin{prop}\label{Proposition4}(Proposition 1 in \cite{Pov99})
Pick $\chi \in (0,1)$ such that
\begin{equation*}
\chi >\max\left(\beta+\alpha_0,1-\left(\frac{\kappa}{d}-\alpha_0\right),
1-\frac{\rho}{d}\right)\label{52}
\end{equation*}
and let
\begin{gather*}
\alpha_1<\min(d(1-\chi),1),\label{53}\\
c_5(a,d,\mu+\nu)>2(1+\gamma(a,d,\mu+\nu)),\label{54}
\end{gather*}
where $\gamma$ and $a$ were introduced in $\eqref{6}$ and $\eqref{7}$, respectively. 
Then we have 
\begin{equation*}
 \begin{split}
\frac{1}{S_t^{\mu,\nu}} {\bold E}_{\epsilon}  \otimes E_0 \bigl[& v(\tau) \,;\,T_{\T}>\tau, \\
& (\mu+\nu)|\mathscr{U}|+ \lambda(\mathscr{U})>c(d,\mu+\nu)+c_5\epsilon^{\alpha_1}\bigr]\\
\le \exp\Bigl\{-(1+\gamma) & t^{\frac{d-\alpha_1}{d+2}}\Bigr\}\label{55}
 \end{split}
\end{equation*}
as $t \to \infty$.
\end{prop}\ \\

\noindent 
If we set 
\begin{equation}
\begin{split}
E =\bigl\{&(\mu+\nu)|\mathscr{U}|+ \lambda(\mathscr{U}) \le c(d,\mu+\nu)+c_5\epsilon^{\alpha_1}, \\
&\lambda_{\omega_1,\omega_2}^{\epsilon}(\T) \le 2c(d,\mu,\nu)\bigr\},\label{69}
\end{split}
\end{equation}
then we have 
\begin{equation*}
\frac{1}{S_t^{\mu,\nu}} {\bold E}_{\epsilon} \otimes E_0 \bigl[ v(\tau) \,;\,T_{\T}>\tau, E^c\bigr]
\le \exp\left\{-(1+\gamma)t^{\frac{d-\alpha_1}{d+2}}\right\}
\end{equation*}
from Proposition \ref{Proposition4}.

The next proposition says that $\mathscr{U}$ is, in a `measurable sense', close to an 
optimal ball of the variational problem in \eqref{4}. 
\begin{prop}\label{Proposition5}(Proposition 2 in \cite{Pov99})
For any $(\omega_1,\omega_2) \in E$, there exists a ball $B$ with radius 
\begin{equation*}
R= \left(\frac{|\mathscr{U}|}{\omega_d}\right)^{\frac{1}{d}}\label{70}
\end{equation*}
such that for large $t$,
\begin{equation}
|\mathscr{U}\setminus B| \le c_6 \epsilon^{\frac{\alpha_1}{20}},
\quad R \le R_0(d,\mu+\nu) + c_6 \epsilon^{\frac{\alpha_1}{2}}\label{71}
\end{equation}
where $R_0$ was defined in $\eqref{5}$ and $c_6(d,\mu+\nu)>0$ is a constant.
\end{prop}\ \\

\noindent 
Thanks to Proposition \ref{Proposition5}, we can introduce for $(\omega_1,\omega_2) \in E$, 
\begin{equation}
B_l: \textrm{the concentric ball to $B$ with radius } R_l=R_0+l\label{85} 
\end{equation}
where $l$ may depends on $t$ and is assumed to satisfy 
\begin{equation*}
l \ge \epsilon^{\alpha_2},\quad \alpha_2 \in \left(0,\alpha_0 \wedge \frac{\alpha_1}{20d}\right).\label{86}
\end{equation*}
It should be pointed out that the ball $B$, and thus $B_l$, corresponding to a 
configuration $(\omega_1,\omega_2) \in E$ need not be unique. 
Since only $B_l\ni 0$ matters in the sequel, we define $B^{(l)}(\omega_1,\omega_2)$ as 
a $B_l$ which contains the origin when $(\omega_1,\omega_2)$ belongs to 
\begin{equation*}
\Omega_1=\{\textrm{There exists a ball $B_l$ with  $0\in B_l$}\} \cap E, \label{87}
\end{equation*}
and otherwise $B^{(l)}(\omega_1,\omega_2)=B(0,R_0+l)$. In particular, we take $l=\epsilon^{\alpha_2}$ 
for Theorem \ref{Theorem1}. Then it follows from \eqref{71} that the radius $R_l$ of $B_l$ satisfies
\begin{equation*}
 \begin{split}
R_0\le R_l &\le R_0+c_6\epsilon^{\frac{\alpha_1}{2}}+\epsilon^{\alpha_2}\\
&\le R_0+(c_6+1)\epsilon^{\alpha_2}\label{88}
 \end{split}
\end{equation*}
and therefore $B^{(l)}(\omega_1,\omega_2)$ has the properties stated in Theorem \ref{Theorem1} with 
\begin{equation*}
\kappa_1=c_6+1,\quad \kappa_2=\alpha_2.\label{89}
\end{equation*}

\subsection{Control of the excursion probability}
In this subsection, we shall derive upper bounds on the probability of 
the excursion of the surviving process from $B$ to $B_l^c$: 
\begin{prop}\label{Proposition6}
There exists a constant $c_7(a,d,\mu+\nu)>0$ such that for $E$ introduced in $\eqref{69}$, 
any $B_l$ defined in $\eqref{85}$ corresponding to $(\omega_1,\omega_2)\in\Omega^2$ 
with $l \le \tau$ and large $t$, we have 
\begin{equation*}
\frac{1}{S_t^{\mu,\nu}} {\bold E}_{\epsilon} \otimes E_0 \bigl[ v(\tau) \,;\,E, T_{\T}>\tau, T_{B_l} \le \tau \bigr]
\le \exp\left\{-c_7l t^{\frac{\alpha_3}{d+2}}\right\} \label{90}
\end{equation*}
where $\alpha_3 = \alpha_0 \wedge (\alpha_1/20d)$.
\end{prop}
This proposition is the (slightly refined) quantitative version of Proposition 3 in \cite{Pov99}. 
The proof will involve estimates on the probability 
of two types of events. These events display two possibilities for surviving 
Brownian motion after it exits $B_l$ : either it returns to $\overline{B}$ 
immediately or stays outside $\overline{B}$ certain amount of time. 
To deal with these events, we shall use two lemmas. 
The first lemma implies that the complement of $B$ almost looks like the `forest set'.
\begin{lem}\label{Lemma1}
Let $(\omega_1,\omega_2) \in E$ and define $c_8=(c_6+1)^{1/d}$, $r_1(\epsilon)=c_8\epsilon^{\alpha_3}$. 
Then there exist constants $c_9(d)>0$, $c_{10}(d)>0$ such that for large $t$,
\begin{gather}
\lambda_{\omega_1,\omega_2}^{\epsilon}(\T \setminus \overline{B}) 
\ge \frac{c_9}{r_1^{2}},\label{92}\\
\sup_{z \in B_l^c}E_z[v(H_{\overline{B}})\,;\,H_{\overline{B}}<T_{\T}]
\le \exp \left\{-c_{10}\frac{l}{r_1} \right\}.\label{93}
\end{gather}
\end{lem}
\noindent{\it Proof. }
This lemma is essentially the same as Lemma 1 in \cite{Pov99} except for dealing with 
longer ``return'' in \eqref{93}. The proof also goes in the same way and we omit the detail. 
\hfill $\square$\\

\eqref{92} and \eqref{93} are related to the estimates on the probability of above events. But they 
are not enough, because we only have rough asymptotics for normalizing constant 
$S_t^{\mu,\nu}$ (cf.\ \eqref{14}). Therefore, we need some cancellation and the next lemma 
meets our need. (This is also essentially the same as Lemma 2 in \cite{Pov99} but we give the proof 
since it seems shorter and self-contained.)
\begin{lem}\label{Lemma2}
\begin{equation*}
\sup_{t\ge 1}\frac{1}{|\T|^2S_t^{\mu,\nu}}{\bold E}_{\epsilon}
\left[\exp\set{-\lambda_{\omega_1,\omega_2}^{\epsilon}(\T)\tau}\right] < \infty.\label{99}
\end{equation*}
\end{lem} 
\noindent{\it Proof. }
We start by introducing the notation
\begin{equation*}
 R_t^{U,V}f(x)=E_x\left[f(Z_t)\exp \left\{-\int_0^t V(Z_s)ds \right\}\,;\,T_U>t \right]\label{101}
\end{equation*}
for a nonempty open set $U$ and a nonnegative function $V$. If $|U|< \infty$ 
and $V$ is locally bounded, then $R_t^{U,V}$ defines a self-adjoint trace class 
semigroup on $L^2(U,dx)$, see for instance (1.3.15) in~\cite{Szn98}. Let 
$\braket{\cdot,\cdot}$ denote the inner product on $L^2(U,dx)$. Using translation 
invariance with respect to $\p$, we find
\begin{equation}
 \begin{split}
S_t^{\mu,\nu}
&={\bold E}_{\epsilon} \otimes E_0 
\left[ \exp \left\{-\int_0^{\tau} V_{\epsilon}(Z_s,\omega_1)ds \right\} \, ; \, 
H_{S_{\epsilon}(\omega_2)}>\tau \right]\\
&=\frac{1}{|\T|}{\bold E}_{\epsilon} \left[ 
\int_{\T} E_x \left[ \exp \left\{-\int_0^{\tau} V_{\epsilon}(Z_s,\omega_1)ds \right\} 
\, ; \, H_{S_{\epsilon}(\omega_2)}>\tau \right]dx \right]\\
&\ge\frac{1}{|\T|}{\bold E}_{\epsilon} \left[ 
\int_{\T\setminus S_{\epsilon}} E_x \bigl[ v(\tau)  \, ; \, T_{\T}>\tau \bigr]dx\,;\,
\lambda_{\omega_1,\omega_2}^{\epsilon}(\T) \le 2c(d,\mu+\nu) \right]\\
&= \frac{1}{|\T|}{\bold E}_{\epsilon} \left[
\braket{R_{\tau}^{\T\setminus S_{\epsilon}(\omega_2), V_{\epsilon}(\cdot,\omega_1)}1,1}\,;\,
\lambda_{\omega_1,\omega_2}^{\epsilon}(\T) \le 2c(d,\mu+\nu)\right]\\
&\ge \frac{1}{|\T|}{\bold E}_{\epsilon} \left[
\braket{\phi,1}^2 \exp\set{-\lambda_{\omega_1,\omega_2}^{\epsilon}(\T)\tau}
\,;\,\lambda_{\omega_1,\omega_2}^{\epsilon}(\T) \le 2c(d,\mu+\nu)\right]\label{102}
 \end{split}
\end{equation}
where $\phi$ is a normalized nonnegative eigenfunction associated with 
$\lambda_{\omega_1,\omega_2}^{\epsilon}(\T)$. Here we have implicitly used 
that $\T$ has finite volume and  $V_{\epsilon}$ is locally bounded for all 
$(\omega_1,\omega_2)\in \Omega^2$. Since $\|R_{s}^{\T\setminus S_{\epsilon}, V_{\epsilon}}\|
_{1 \to \infty}\le (2\pi s)^{-\frac{d}{2}}$ for all $s>0$, we have
\begin{equation*}
 \begin{split}
0 &\le \exp\set{-\lambda_{\omega_1,\omega_2}^{\epsilon}(\T)} \phi(z)\\ 
&= R_{1}^{\T\setminus S_{\epsilon}(\omega_2),V_{\epsilon}(\cdot,\omega_1)}\phi(z)\\
&\le (2\pi)^{-\frac{d}{2}}\braket{\phi,1}\label{103}
 \end{split}
\end{equation*}
and therefore for $(\omega_1,\omega_2)\in \{\lambda_{\omega_1,\omega_2}^{\epsilon}(\T) \le 2c(d,\mu+\nu)\}$,
\begin{equation*}
 \begin{split}
\braket{\phi,1}^2&\ge (2\pi)^d\exp\set{-4c(d,\mu+\nu)}\sup_{z\in\R^d}\phi(z)^2\\
&\ge (2\pi)^d\exp\set{-4c(d,\mu+\nu)}|\T|^{-1}\label{104}
 \end{split}
\end{equation*}
where we have used $\mathrm{supp}\,\phi\subset\T$ together with $\|\phi\|_2=1$ 
in the second inequality. Coming back to \eqref{102}, we have shown that for $t \ge 1$,
\begin{equation*}
 \begin{split}
S_t^{\mu,\nu}\ge \frac{c_{11}(d,\mu+\nu)}{|\T|^2}{\bold E}_{\epsilon} 
\Bigl[&\exp\set{-\lambda_{\omega_1,\omega_2}^{\epsilon}(\T)\tau} \,; \\
&\lambda_{\omega_1,\omega_2}^{\epsilon}(\T) \le 2c(d,\mu+\nu)\Bigr].\label{105}
 \end{split}
\end{equation*}
We can drop $\{\lambda_{\omega_1,\omega_2}^{\epsilon}(\T) \le 2c(d,\mu+\nu)\}$ since 
\begin{equation*}
\begin{split}
&\frac{1}{S_t^{\mu,\nu}}{\bold E}_{\epsilon} \Bigl[\exp
\set{-\lambda_{\omega_1,\omega_2}^{\epsilon}(\T)\tau}\,;\,\lambda_{\omega_1,\omega_2}^{\epsilon}(\T)>2c(d,\mu+\nu)\Bigr]\\ 
\le \, &\frac{1}{S_t^{\mu,\nu}}\exp\set{-2c(d,\mu+\nu)\tau} \to 0
\end{split}
\end{equation*}
as $t \to \infty$ from \eqref{14}. Thus the claim of Lemma 2 follows.
\hfill $\square$\\

We are now ready to prove Proposition \ref{Proposition6}.\\

\noindent{\it{Proof of Proposition \ref{Proposition6}.}} 
Let $\theta_{\cdot}$ denote the canonical shift and define 
\begin{equation*}
 \begin{split}
\set{T_{B_l}\le\tau<T_{\T}}\subset 
&\,\set{2T_{B_l}\le\tau<T_{\T},H_{\overline{B}}\circ\theta_{T_{B_l}}>lr_1}\\
&\,\cup\set{2T_{B_l}\le\tau<T_{\T},H_{\overline{B}}\circ\theta_{T_{B_l}}\le lr_1}\\
&\,\cup\set{\frac{\tau}{2}<T_{B_l}\le\tau,T_{\T}>\tau}\\
=&\, A_1\cup A_2\cup A_3\label{115}
 \end{split}
\end{equation*}
for $t \ge 1$ and $(\omega_1,\omega_2) \in E$. 
Here note that our choice $l\le \tau$ implies $\tau/2+lr_1<\tau$ for large $t$. 
We start with $A_1$. Pick $(\omega_1,\omega_2) \in E$ and write 
\begin{equation}
 \begin{split}
E_0 \left[v(\tau)\,;\, A_1\right]
&=\,\sum_{k=0}^{[\tau/2]}E_0 \bigl[v(\tau)\,;\, k \le T_{B_l} <k+1 , A_1 \bigr]\\
&=\,\sum_{k=0}^{[\tau/2]}E_0 \bigl[v(\tau)\,;\, E_k , A_1 \bigr].\label{116}
 \end{split}
\end{equation}
Applying strong Markov property at $T_{B_l}+lr_1$ and $T_{B_l}$ and using (3.1.9) of \cite{Szn98}, 
we find for $1\le k \le [\tau/2]$ and large $t$ that
\begin{equation}
 \begin{split}
&E_0 \left[v(\tau)\,;\, E_k,A_1 \right]\\
\le\, & c(d)\left(1+\left(\lambda_{\omega_1,\omega_2}^{\epsilon}(\T)\tau\right)^{\frac{d}{2}}\right)
\exp\set{-\lambda_{\omega_1,\omega_2}^{\epsilon}(\T)(\tau-(k+1+lr_1))}\\
&E_0\left[v(T_{B_l}+lr_1)\,;\,E_k,H_{\overline{B}}\circ\theta_{T_{B_l}}>lr_1,
T_{\T}>T_{B_l},T_{\T}\circ\theta_{T_{B_l}}>lr_1\right]\\
\le\, & c(d)\left(1+\left(2c(d,\mu+\nu)\tau\right)^{\frac{d}{2}}\right)
\exp\set{-\lambda_{\omega_1,\omega_2}^{\epsilon}(\T)(\tau-(k+1+lr_1))}\\
&\sup_{z\in B^c}E_z\bigl[v(lr_1)\,;\,H_{\overline{B}}>lr_1,T_{\T}>lr_1\bigr]
E_0\bigl[v(T_{B_l})\,;\,E_k,T_{\T}>T_{B_l}\bigr]\\
\le\, & c(d)^2\left(1+\left(2c(d,\mu+\nu)\tau\right)^{\frac{d}{2}}\right)
\exp\set{-\lambda_{\omega_1,\omega_2}^{\epsilon}(\T)(\tau-(k+1+lr_1))}\\
&\left(1+\left(\lambda_{\omega_1,\omega_2}^{\epsilon}(\T\setminus\overline{B})lr_1\right)^{\frac{d}{2}}\right)
\exp\set{-\lambda_{\omega_1,\omega_2}^{\epsilon}(\T \setminus \overline{B})lr_1}E_0\bigl[v(k)\,;\,T_{\T}>k\bigr]\\
\le\, & c_{12}\exp\left\{-\lambda_{\omega_1,\omega_2}^{\epsilon}(\T)(\tau-lr_1)
-\frac{1}{2}\lambda_{\omega_1,\omega_2}^{\epsilon}(\T \setminus \overline{B})lr_1\right\}\label{117}
 \end{split}
\end{equation}
for some constant $c_{12}(d,\mu+\nu)>0$. Here we have used $\lambda_{\omega_1,\omega_2}^{\epsilon}(\T)\le
2c(d,\mu+\nu)$ for $(\omega_1,\omega_2)\in E$ and 
\begin{equation}
\lambda_{\omega_1,\omega_2}^{\epsilon}(\T \setminus \overline{B})lr_1
\ge \frac{c_9l}{r_1} \ge \frac{c_9}{c_8}t^{\frac{\alpha_3-\alpha_2}{d+2}}\label{117'}
\end{equation}
from \eqref{92} and our choice $l \ge t^{-\alpha_2/(d+2)}$. Coming back to \eqref{116} and using \eqref{117'} again, we have
\begin{equation*}
 \begin{split}
&\frac{1}{S_t^{\mu,\nu}}{\bold E}_{\epsilon} \otimes E_0 \bigl[v(\tau) \,;\, E \cap A_1\bigr]\\
\le\, & c_{12}\tau^{d+1}|\T|^2 \exp \left\{-\frac{c_9l}{2r_1}+(c(d,\mu+\nu)+1)lr_1\right\}\\
& \sup_{t\ge 1}\frac{1}{|\T|^2S_t^{\mu,\nu}}{\bold E}_{\epsilon}
\left[\exp\set{-\lambda_{\omega_1,\omega_2}^{\epsilon}(\T)\tau}\,;\,E \right].\label{118}
 \end{split}
\end{equation*}
Therefore for large $t$, Lemma \ref{Lemma2} gives us 
\begin{equation*}
\frac{1}{S_t^{\mu,\nu}}{\bold E}_{\epsilon} \otimes 
E_0 \bigl[v(\tau)\,;\,E \cap A_1\bigr]\le \exp\left\{-c_9\frac{l}{r_1}\right\}\label{119}
\end{equation*}
with slightly smaller $c_9$. 

Next, we shall deal with $A_2$. As in \eqref{116} we write 
\begin{equation}
E_0 \bigl[v(\tau)\,;\, A_2\bigr]=\sum_{k=0}^{[\tau/2]}E_0 \bigl[v(\tau)\,;\, E_k , A_2 \bigr].\label{120}
\end{equation}
for $(\omega_1,\omega_2)\in E$. Then we have for $0\le k \le[\tau/2]$ and 
large enough $t$, 
\begin{equation}
 \begin{split}
&E_0 \left[v(\tau)\,;\, E_k,A_2 \right]\\
\le\, & c(d)\left(1+\left(\lambda_{\omega_1,\omega_2}^{\epsilon}(\T)\tau)^{\frac{d}{2}}\right)\right)
\exp\set{-\lambda_{\omega_1,\omega_2}^{\epsilon}(\T)(\tau-(k+1+lr_1))}\\
&E_0\left[v(T_{B_l}+lr_1)\,;\,E_k,H_{\overline{B}}
\circ\theta_{T_{B_l}}\le lr_1,T_{\T}>T_{B_l},T_{\T}\circ\theta_{T_{B_l}}>lr_1\right]\\
\le\, & c(d)\left(1+\left(2c(d,\mu+\nu)\tau)^{\frac{d}{2}}\right)\right)
\exp\set{-\lambda_{\omega_1,\omega_2}^{\epsilon}(\T)(\tau-(k+1+lr_1))}\\
&\sup_{z\in B^c}E_z\bigl[v(H_{\overline{B}})\,;\,H_{\overline{B}}\le lr_1<T_{\T}\bigr]
E_0\bigl[v(T_{B_l})\,;\,E_k,T_{\T}>T_{B_l}\bigr]\\
\le\, & c_{13}\tau^d\sup_{z\in B^c}E_z\bigl[v(H_{\overline{B}})\,;\,H_{\overline{B}}<T_{\T}\bigr]
\exp\left\{-\lambda_{\omega_1,\omega_2}^{\epsilon}(\T)(\tau-lr_1)\right\}\label{121}
 \end{split}
\end{equation}
for some constant $c_{13}(d,\mu+\nu)>0$. Coming back to \eqref{120} and using \eqref{93}, 
we find
\begin{equation*}
 \begin{split}
&\frac{1}{S_t^{\mu,\nu}}{\bold E}_{\epsilon} \otimes E_0 [v(\tau) \,;\, E \cap A_2]\\
\le\, & c_{13}\tau^{d+1}|\T|^2 \exp \left\{-\frac{c_{10}l}{r_1}+(c(d,\mu+\nu)+1)lr_1\right\}\\
&\sup_{t\ge 1}\frac{1}{|\T|^2S_t^{\mu,\nu}}{\bold E}_{\epsilon}
\left[\exp\set{-\lambda_{\omega_1,\omega_2}^{\epsilon}(\T)\tau}\,;\,E \right].\label{122}
 \end{split}
\end{equation*}
Therefore it follows as before that 
\begin{equation*}
\frac{1}{S_t^{\mu,\nu}}{\bold E}_{\epsilon} \otimes 
E_0 [v(\tau)\,;\,E \cap A_2]\le \exp\left\{-c_{10}\frac{l}{r_1}\right\}\label{123}
\end{equation*}
with slightly smaller $c_{10}$.

As for $A_3$, observe that on $\{\tau/2 < T_{B_l} \le \tau\}$ the reversed path 
starting at $Z_{\tau}$ exits $B_l$ before time $\tau/2$. Since the estimates \eqref{117} and \eqref{121} 
do not depend on the starting point 0, 
we see that $E_0 [v(\tau)\,;\, A_3]$ is bounded above by $[\tau/2]$ times the sum of the 
right hand side in \eqref{117} and \eqref{121}, respectively. So, we have the same upper bound on 
${\bold E}_{\epsilon} \otimes E_0 [v(\tau)\,;\,
E \cap A_3]$ as on $A_1 \cup A_2$. The proof of Proposition \ref{Proposition6} is now complete.
\hfill $\square$

\subsection{Proof of the upper estimates}
Now we are ready to prove Theorem \ref{Theorem1} and the upper estimates of Theorem \ref{Theorem2} and \ref{Theorem3}. 
To this end, we are going to give an upper bound on the probability that (scaled) surviving 
process leaves $B^{(l)}(\omega_1,\omega_2)$ before time $\tau$. 
We first note that for $(\omega_1,\omega_2)\in\Omega_1^c\cap E$ the starting 
point of the process is not contained in any $B_l$ and thus $T_{B_l}=0$. 
Therefore, no matter which $B_{l}$ we pick for $B^{(l)}$ on $\Omega_1^c\cap E$, we have 
\begin{equation*}
 \begin{split}
\set{T_{B^{(l)}} \le \tau} &\subset
\set{T_{B^{(l)}}\le\tau,T_{\T}>\tau}\cup\set{T_{\T} \le \tau}\\
&\subset\set{T_{B^{(l)}}\le\tau, T_{\T}>\tau,E}\cup\set{T_{\T}>\tau,E^c}\cup
\set{T_{\T} \le \tau}.\label{124}
 \end{split}
\end{equation*}
Consequently, we find for large $t$ that
\begin{equation}
 \begin{split}
&Q_t^{\mu,\nu}\left(T_{t^{1/(d+2)}B^{(l)}(\omega_1,\omega_2)}\le t\right)\\
\le &\,\frac{1}{S_t^{\mu,\nu}} {\bold E}_{\epsilon} 
\otimes E_0\bigl[v(\tau)\,;\,T_{B^{(l)}}\le\tau, T_{\T}>\tau, E \bigr]\\
&+\frac{1}{S_t^{\mu,\nu}} {\bold E}_{\epsilon} 
\otimes E_0\bigl[v(\tau)\,;\,E^c \bigr]\\
&+\frac{1}{S_t^{\mu,\nu}} {\bold E}_{\epsilon} 
\otimes E_0\bigl[v(\tau)\,;\,T_{\T} \le \tau \bigr]\\
\le &\,\exp\left\{-c_{14}\Bigl(lt^{\frac{\alpha_3}{d+2}}\wedge t^{\frac{d-\alpha_1}{d+2}}\Bigr) \right\}\label{125}
 \end{split}
\end{equation}
with $c_{14}=c_7 \wedge (1+\gamma)$ using \eqref{50}, Proposition \ref{Proposition4} 
and Proposition \ref{Proposition6}. \\

\noindent\emph{Proof of Theorem \ref{Theorem1}.} We set $l=\epsilon^{\alpha_2}$ 
as previously stated at the end of section 3.3. Then we have 
\begin{equation*}
lt^{\frac{\alpha_3}{d+2}}\wedge t^{\frac{d-\alpha_1}{d+2}}=
t^{\frac{\alpha_3-\alpha_2}{d+2}}\wedge t^{\frac{d-\alpha_1}{d+2}} 
\end{equation*}
and therefore the right hand side of \eqref{125} converges to $0$ as $t\to\infty$, 
which proves Theorem \ref{Theorem1}. \hfill $\square$ \\

\noindent\emph{Proof of the upper estimate of $\eqref{9}$}. 
Since on $\{T_{t^{1/(d+2)}B^{(l)}(\omega_1,\omega_2)}>t\}$ we have 
\begin{equation}
W_t^C \subset a\textrm{-neighborhood of }t^{\frac{1}{d+2}}B^{(l)}(\omega_1,\omega_2)\label{126}
\end{equation}
and the volume of the right hand side of \eqref{126} is smaller than
\begin{equation*}
t^{\frac{d}{d+2}}\omega_d \left(R_0+\kappa_1t^{-\frac{\kappa_2}{d+2}}+at^{-\frac{1}{d+2}}\right)^d
\sim t^{\frac{d}{d+2}}\omega_d R_0^d \quad (t \to \infty),\label{127}
\end{equation*}
the upper estimate of Theorem \ref{Theorem2} follows. \hfill $\square$\\ 

\noindent\emph{Proof of Theorem \ref{Theorem3}.} It is enough to show that for arbitrary $\eta >0$ 
and $l > 2R_0$,
\begin{equation}
Q_t^{\mu,\nu}\left(T_{t^{1/(d+2)}B(0,l)}\le t\right) \le \exp \set{-\eta l}\label{127'}
\end{equation}
when $t$ is large enough. 
First of all, we can find a constant $M>0$ such that for all $l > Mt^{d/(d+2)}$,
\begin{equation*}
 P_0\left(T_{t^{1/(d+2)}B(0,l)}\le t\right) \le \exp \set{-\eta l -2c(d,\mu+\nu)t^{\frac{d}{d+2}}}
\end{equation*}
using a standard Brownian estimate. From this and \eqref{14}, we have \eqref{127'} in this case.
Next, for $2R_0 < l \le M^2 t^{(d-1)/(d+2)}$ we can use \eqref{125} to derive
\begin{equation*}
 \begin{split}
  Q_t^{\mu,\nu}\left(T_{t^{1/(d+2)}B(0,l)}\le t\right) 
  &\le Q_t^{\mu,\nu}\left(T_{t^{1/(d+2)}B^{(l/2-R_0)}(\omega_1,\omega_2)}\le t\right)\\
  &\le \exp\left\{-c_{14}\left(\left(\frac{l}{2}-R_0\right)t^{\frac{\alpha_3}{d+2}}\wedge 
  \frac{l}{M^2} t^{\frac{1-\alpha_1}{d+2}} \right) \right\}.\label{127''}
 \end{split}
\end{equation*}
This implies \eqref{127'} when $t$ is large. 
In the remaining case $M^2 t^{(d-1)/(d+2)} < l \le Mt^{d/(d+2)}$, we shall use large deviation estimate 
in \cite{Szn95a} again to show
\begin{equation}
 \limsup_{t \to \infty}t^{-\frac{d}{d+2}}\log Q_t^{\mu,\nu}\left(\sup_{0 \le s \le t}|Z_s| > xt^{\frac{d}{d+2}}\right)
 \le - \inf_{y \not\in B(0,x)}\beta_0(y). \label{128}
\end{equation}
Strictly speaking, \cite{Szn95a} deals with large deviation estimates for $t^{-d/(d+2)}Z_t$ but the proof is obviously
 applicable to above version. 
The only property of $\beta_0$ we need here is that it is a norm. 
We can deduce from \eqref{128} that for large $t$,
\begin{equation}
 \begin{split}
  Q_t^{\mu,\nu}\left(T_{t^{1/(d+2)}B(0,l)}\le t \right) 
  &\le Q_t^{\mu,\nu} \left(\sup_{0 \le s \le t}|Z_s| > M^2 t^{\frac{d}{d+2}} \right)\\
  &\le \exp\left\{-\frac{1}{2} M^2 t^{\frac{d}{d+2}} \inf_{y \not\in B(0,1)}\beta_0(y) \right\}\\
  &\le \exp\left\{-\frac{1}{2} M l \inf_{y \not\in B(0,1)}\beta_0(y) \right\}.\label{129}
 \end{split}
\end{equation}
This implies \eqref{127'} making $M$ larger if necessary. 

Now we have \eqref{127'} for all $l > 2R_0$ and the proof of Theorem 3 is completed. \hfill $\square$ 

\section*{Acknowledgement}
The auther would like to thank Professor Nobuo Yoshida for suggesting this very interesting problem, 
a lot of  helpful discussions, and careful reading of the early version of the manuscript.  
He is also grateful to the referee for many constructive suggestions.


\begin{thebibliography}{1}

\bibitem{Ber03}
Marcel Berger.
\newblock {\em A panoramic view of {R}iemannian geometry}.
\newblock Springer-Verlag, Berlin, 2003.

\bibitem{DZ98}
Amir Dembo and Ofer Zeitouni.
\newblock {\em Large deviations techniques and applications}, volume~38 of {\em
  Applications of Mathematics (New York)}.
\newblock Springer-Verlag, New York, second edition, 1998.

\bibitem{DV75}
M.~D. Donsker and S.~R.~S. Varadhan.
\newblock Asymptotics for the {W}iener sausage.
\newblock {\em Comm. Pure Appl. Math.}, 28(4):525--565, 1975.

\bibitem{Pov99}
Tobias Povel.
\newblock Confinement of {B}rownian motion among {P}oissonian obstacles in
  {${\bf R}\sp d,\ d\ge3$}.
\newblock {\em Probab. Theory Related Fields}, 114(2):177--205, 1999.

\bibitem{Szn91b}
Alain-Sol Sznitman.
\newblock On the confinement property of two-dimensional {B}rownian motion
  among {P}oissonian obstacles.
\newblock {\em Comm. Pure Appl. Math.}, 44(8-9):1137--1170, 1991.

\bibitem{Szn95a}
Alain-Sol Sznitman.
\newblock Annealed {L}yapounov exponents and large deviations in a {P}oissonian
  potential. {I}, {II}.
\newblock {\em Ann. Sci. \'Ecole Norm. Sup. (4)}, 28(3):345--370, 371--390,
  1995.

\bibitem{Szn98}
Alain-Sol Sznitman.
\newblock {\em Brownian motion, obstacles and random media}.
\newblock Springer Monographs in Mathematics. Springer-Verlag, Berlin, 1998.

\end{thebibliography}
\end{document}